# Direction curves of tangent indicatrix of a curve

**Burak Şahiner**

*Manisa Celal Bayar University, Faculty of Arts and Sciences, Department of Mathematics, 45140 Manisa, Turkey.*
*E-mail: burak.sahiner@cbu.edu.tr*

**Abstract**
In this paper, we define some new associated curves as integral curves of a vector field generated by Frenet vectors of tangent indicatrix of a curve in Euclidean 3-space. We give some relationships between curvatures of these curves. By using these associated curves, we give some methods to construct helices and slant helices from some special spherical curves such as circles on unit sphere, spherical helices, and spherical slant helices. Finally, we give some related examples.

## 1. Introduction

The study of curves has wide range applications in many areas. While helix curves arise in DNA double, nanosprings, and carbon nanotubes, Bezier curves are used in computer aided geometric design. The computation of curves on manifolds also finds applications in some other fields (for details, see (Eisenhart 1940; Do Carmo, 1976; Struik, 1988; Allasia et al., 2018)). One of the interesting topics in the theory of curves is associated curves. Two or more curves are called associated curves if there is some mathematical relationship between them. Despite being studied for many years, the topic of associated curves is still one of the most interesting research areas in the theory of curves. Some important characterizations and some intrinsic properties such as Frenet vectors and curvatures of a curve can be found by using its associated curve. Involute-evolute curve couples, Bertrand curves, and Mannheim partner curves are some well-known associated curves.

Recently, Choi and Kim (2012) defined a new type of associated curves called direction curve as an integral curve of a vector field generated by Frenet vectors of a given curve. Then, they found some relations between curvatures of these associated curves. Moreover, they gave a canonical method to construct general helices and slant helices by using these curves. After the paper of Choi and Kim (2012), many authors have taken interest in the study of this new type of associated curve. Choi et al. (2012), and Qian and Kim (2015) studied direction curves in Minkowski space $E_1^3$. Körpınar et al. (2013) studied direction curves by using Bishop frame instead of Frenet frame in their work. Macit and Düldül (2014) defined $W$-direction curve by using unit Darboux vector field $W$ of a given curve and introduced $V$-direction curve of a given curve on a surface by using Darboux frame. They also examined direction curve of a given curve in four-dimensional Euclidean space. Kızıltuğ and Önder (2015) gave general definition of direction curves in a three-dimensional compact Lie group.

In this paper, we first define $X$-direction curve as an integral curve of unit vector field $X$ generated by Frenet vectors of tangent indicatrix of a given curve. By using this definition, we define evolute-direction curves, Bertrand-direction curves, and Mannheim-direction curves of tangent indicatrix of the given curve. For each one of all these direction curves, we find some relations between Frenet vectors and curvatures. Moreover, we give some methods to construct helices and slant helices from some special spherical curves which are circles on unit sphere, spherical helices, and spherical slant helices. Finally, we give some examples of these associated curves.

## 2. Preliminaries

In this section, we give some basic concepts of a curve in Euclidean 3-space, tangent indicatrix, general helix, slant helix, and some associated curves.

Let $\alpha : I \subset IR \to IR^3$ be a unit speed curve in Euclidean 3-space. Three orthonormal vectors of Frenet frame $\{T, N, B\}$ along the curve $\alpha$ can be defined as follows

$$T = \frac{d\alpha}{ds},\ N = (1/\kappa)T',\ B = T \times N,$$

where $T$ is unit tangent vector field, $N$ is unit principal normal vector field, $B$ is unit binormal vector field, $\kappa = \|T'\|$ is called curvature which measures the amount of deviation of the curve from the tangent line, $s$ is the arc-length parameter of $\alpha$. The Frenet derivative formulas can also be given in matrix form as

$$\frac{d}{ds}\begin{bmatrix} T \\ N \\ B \end{bmatrix} = \begin{bmatrix} 0 & \kappa & 0 \\ -\kappa & 0 & \tau \\ 0 & -\tau & 0 \end{bmatrix}\begin{bmatrix} T \\ N \\ B \end{bmatrix},$$

where $\tau = -\langle B', N \rangle$ is called torsion which measures the amount of deviation of the curve from the osculating plane spanned by the vectors $T$ and $N$ (Do Carmo, 1976; O'Neill, 2006).

A curve $\alpha$ is called helix or general helix if its unit tangent vector field $T(s)$ makes a constant angle with a fixed straight line. $\alpha$ is a general helix if and only if

$$f(s) = \frac{\tau}{\kappa}(s) \tag{2.1}$$

is a constant function (Struik, 1988).

Similarly, a curve $\alpha$ is called slant helix if its unit principal normal vector field $N(s)$ makes a constant angle with a fixed straight line. $\alpha$ is a slant helix if and only if

$$\sigma(s) = \frac{\kappa^2}{(\kappa^2 + \tau^2)^{3/2}}\left(\frac{\tau}{\kappa}\right)'(s) \tag{2.2}$$

is a constant function (Izumiya and Takeuchi, 2002).

Let $\alpha : I \subset IR \to IR^3$ and $\beta : I \subset IR \to IR^3$ be two curves in Euclidean 3-space with Frenet frames $\{T_\alpha, N_\alpha, B_\alpha\}$ and $\{T_\beta, N_\beta, B_\beta\}$, respectively. $\beta$ is called evolute (resp., Bertrand mate, Mannheim partner) curve of $\alpha$ if and only if $N_\beta = T_\alpha$ (resp., $N_\beta = N_\alpha$, $N_\beta = B_\alpha$) (Eisenhart, 1940; Wang and Liu, 2007).

Let $\alpha : I \subset IR \to IR^3$ be unit speed curve in Euclidean 3-space with Frenet frame apparatus $\{T, N, B, \kappa, \tau\}$. The unit tangent vectors along the curve $\alpha(s)$ generate a curve $\alpha_T = T$ on the unit sphere about the origin. The curve $\alpha_T$ is called tangent indicatrix of the curve $\alpha$. If the Frenet frame apparatus of $\alpha_T$ is $\{T_T, N_T, B_T, \kappa_T, \tau_T\}$, then the Frenet derivative formulas of $\alpha_T$ can be given as follows

$$\frac{d}{ds_T}\begin{bmatrix} T_T \\ N_T \\ B_T \end{bmatrix} = \begin{bmatrix} 0 & \kappa_T & 0 \\ -\kappa_T & 0 & \tau_T \\ 0 & -\tau_T & 0 \end{bmatrix}\begin{bmatrix} T_T \\ N_T \\ B_T \end{bmatrix}, \tag{2.3}$$

where

$$T_T = N,\ N_T = \frac{-T + fB}{\sqrt{1+f^2}},\ B_T = \frac{fT + B}{\sqrt{1+f^2}},$$

and

$$s_T = \int \kappa(s)ds\ ,\ \kappa_T = \sqrt{1+f^2}\ ,\ \ \tau_T = \sigma\sqrt{1+f^2}\ ,$$

where $s_T$ is natural representation of the tangent indicatrix $\alpha_T$ of the curve $\alpha$, $f$ and $\sigma$ are functions given by equations (2.1) and (2.2), respectively (Ali, 2012).

## 3. Direction curves of tangent indicatrix

Let $\alpha : I \subset IR \to E^3$ be a curve in Euclidean 3-space with Frenet apparatus $\{T, N, B, \kappa, \tau\}$ and $\alpha_T$ be tangent indicatrix of $\alpha$ with Frenet apparatus $\{T_T, N_T, B_T, \kappa_T, \tau_T\}$. Consider a vector field $X$ expressed by

$$X(s_T) = x(s_T)T_T(s_T) + y(s_T)N_T(s_T) + z(s_T)B_T(s_T)\ , \tag{3.1}$$

where $x, y$ and $z$ are real functions of $s_T$ which is natural representation of the tangent indicatrix of the curve $\alpha$. If we assume that the vector field $X$ is unit, we get

$$x^2(s_T) + y^2(s_T) + z^2(s_T) = 1. \tag{3.2}$$

By differentiating equation (3.2), we have

$$x(s_T)x'(s_T) + y(s_T)y'(s_T) + z(s_T)z'(s_T) = 0\ . \tag{3.3}$$

Now we can give the definition of $X$ -direction curve of tangent indicatrix of a curve as follows.

**Definition 3.1.** Let $\alpha$ be a curve in Euclidean 3-space, $\alpha_T$ be tangent indicatrix of the curve $\alpha$, and $X$ be a unit vector field satisfying the equations (3.1) and (3.2). The integral curve $\beta : I \subset IR \to IR^3$ of $X$ is called $X$ -direction curve of $\alpha_T$ or $X$ -direction curve of tangent indicatrix of $\alpha$. $\alpha_T$ is also called $X$ -donor curve of $\beta$.

Let $\beta$ be $X$ -direction curve of the tangent indicatrix $\alpha_T$ of the curve $\alpha$ with Frenet frame $\{T_\beta, N_\beta, B_\beta\}$ and curvatures $\{\kappa_\beta, \tau_\beta\}$. Since the vector field $X$ is unit and $\frac{d\beta}{ds_\beta} = T_\beta = X$, the arc-length parameter $s_\beta$ of $\beta$ equals to $s_T + c$, where $c$ is a constant real number. Without loss of generality, we assume that $s_\beta = s_T$. By differentiating equation (3.1) and using the Frenet formulae in Equation (2.3), we have

$$\kappa_\beta N_\beta = (x' - y\kappa_T)T_T + (y' + x\kappa_T - z\tau_T)N_T + (z' + y\tau_T)B_T\ . \tag{3.4}$$

By using equality (3.4), we can give definitions of evolute-direction, Bertrand-direction, and Mannheim-direction curves of the tangent indicatrix $\alpha_T$, and study some properties of these curves.

## 4. Evolute-direction curves of tangent indicatrix

In this section, we define evolute-direction curves of tangent indicatrix of a curve and obtain some relations between these curves.

**Definition 4.1.** Let $\alpha$ be a curve in Euclidean 3-space, $\alpha_T$ be tangent indicatrix of $\alpha$, and $\beta$ be $X$-direction curve of $\alpha_T$. If $\beta$ is an evolute of $\alpha_T$, then $\beta$ is called evolute-direction curve of $\alpha_T$. $\alpha_T$ is also called involute-donor curve of $\beta$.

**Proposition 4.2.** Let $\alpha$ be a curve in Euclidean 3-space and $\alpha_T$ be tangent indicatrix of $\alpha$. $\beta$ is evolute-direction curve of $\alpha_T$ if and only if the functions in (3.1) are as follows

$$x(s_T)=0,\quad y(s_T)=\sin\left(\int \tau_T(s_T)ds_T\right),\ z(s_T)=\cos\left(\int \tau_T(s_T)ds_T\right).$$

**Proof.** From the definition of evolute curves, we know that $N_\beta = T_T$. By using equality (3.4), we have the following system of differential equations

$$\left.\begin{aligned} x' - y\kappa_T &= \kappa_\beta \\ y' + x\kappa_T - z\tau_T &= 0 \\ z' + y\tau_T &= 0 \end{aligned}\right\} \tag{4.1}$$

Multiplying the first, second and third equations in (4.1) with $x$, $y$ and $z$, respectively, adding the results, and using equation (3.3), we obtain $x=0$. By substituting $x=0$ into system (4.1), the solution is found as follows

$$\left\{ x(s_T)=0,\ y(s_T)=\sin\left(\int \tau_T(s_T)ds_T\right),\ z(s_T)=\cos\left(\int \tau_T(s_T)ds_T\right) \right\}.$$

From Proposition 4.2, we have a method to construct a unit speed evolute curve from a unit spherical curve. This construction can be achieved just by using the Frenet vectors $N_T, B_T$ and the torsion $\tau_T$ of a unit spherical curve $\alpha_T$ which is tangent indicatrix of a curve $\alpha$.

Now we can give the following relationships between curvatures of the tangent indicatrix $\alpha_T$ of the curve $\alpha$ and its evolute-direction curve $\beta$.

**Theorem 4.3.** Let $\alpha$ be a curve in Euclidean 3-space, $\alpha_T$ be tangent indicatrix of $\alpha$, and $\beta$ be evolute-direction curve of $\alpha_T$. The relations between curvatures of the curves $\alpha_T$ and $\beta$ can be given as follows

$$\kappa_\beta = -\kappa_T \sin\left(\int \tau_T(s_T)ds_T\right),\quad \tau_\beta = \kappa_T \cos\left(\int \tau_T(s_T)ds_T\right)$$

and

$$\kappa_T = \sqrt{\kappa_\beta^{\,2} + \tau_\beta^{\,2}},\quad \tau_T = \frac{\kappa_\beta^{\,2}}{\kappa_\beta^{\,2} + \tau_\beta^{\,2}}\left(\frac{\tau_\beta}{\kappa_\beta}\right)'.$$

**Proof.** From the first equation of system (4.1), we can see that $\kappa_\beta = -\kappa_T \sin\left(\int \tau_T(s_T)ds_T\right)$. From Proposition (4.2), we find $X = T_\beta = \sin\left(\int \tau_T(s_T)ds_T\right)N_T + \cos\left(\int \tau_T(s_T)ds_T\right)B_T$, and from the definition of evolute curves, we have $N_\beta = T_T$. The unit binormal vector $B_\beta$ of the evolute-direction curve $\beta$ can be found as $B_\beta = T_\beta \times N_\beta = \cos\left(\int \tau_T(s_T)ds_T\right)N_T - \sin\left(\int \tau_T(s_T)ds_T\right)B_T$. By differentiating unit binormal vector $B_\beta$ and using the third equality of (2.3), we get $\tau_\beta = \kappa_T \cos\left(\int \tau_T(s_T)ds_T\right)$.

On the other hand, by differentiating the equality $N_\beta = T_T$ and using (2.3), we get $\kappa_T = \sqrt{\kappa_\beta{}^2 + \tau_\beta{}^2}$. By using the curvatures of the curve $\beta$, we have $\frac{\tau_\beta}{\kappa_\beta} = -\cot\left(\int \tau_T ds_T\right)$, and so $\int \tau_T ds_T = -\operatorname{arc}\cot\left(\frac{\tau_\beta}{\kappa_\beta}\right)$. By differentiating the last equality, we get $\tau_T = \frac{\kappa_\beta{}^2}{\kappa_\beta{}^2 + \tau_\beta{}^2}\left(\frac{\tau_\beta}{\kappa_\beta}\right)'$.

□

From Theorem 4.3, we have the following corollary.

**Corollary 4.4.** Let $\alpha$ be a curve in Euclidean 3-space, $\alpha_T$ be tangent indicatrix of $\alpha$, and $\beta$ be evolute-direction curve of $\alpha_T$. The following relations between curvatures of the curves $\alpha_T$ and $\beta$ are satisfied

$$\frac{\tau_\beta}{\kappa_\beta} = \cot\left(\int \tau_T ds_T\right) \quad \text{and} \quad \frac{\tau_T}{\kappa_T} = \frac{\kappa_\beta{}^2}{\left(\kappa_\beta{}^2 + \tau_\beta{}^2\right)^{3/2}}\left(\frac{\tau_\beta}{\kappa_\beta}\right)'.$$

From Corollary 4.4, we can give the following theorem without any proof.

**Theorem 4.5.** Let $\alpha$ be a curve in Euclidean 3-space, $\alpha_T$ be tangent indicatrix of $\alpha$, and $\beta$ be evolute-direction curve of $\alpha_T$. Then

**i)** $\alpha_T$ is a circle or a part of circle on unit sphere if and only if $\beta$ is a helix.

**ii)** $\alpha_T$ is a spherical helix if and only if $\beta$ is a slant helix.

Theorem 4.5 gives a method to construct a helix and a slant helix from a circle or a part of circle on unit sphere and a spherical helix, respectively.

Now, by using relations between curvatures of the curve $\alpha$ and its tangent indicatrix $\alpha_T$, we can give the following relations between curvatures of the curve $\alpha$ and evolute-direction curve $\beta$ of $\alpha_T$.

**Corollary 4.6.** Let $\alpha$ be a curve in Euclidean 3-space with nonzero curvature, $\alpha_T$ be tangent indicatrix of $\alpha$, and $\beta$ be evolute-direction curve of $\alpha_T$. The following relations between curvatures of the curves $\alpha$ and $\beta$ are satisfied

$$\frac{\tau_\beta}{\kappa_\beta} = \cot\left(\int \frac{\kappa^3}{\left(\kappa^2+\tau^2\right)^{3/2}}\left(\frac{\tau}{\kappa}\right)'\sqrt{1+\left(\frac{\tau}{\kappa}\right)^2}\,ds\right) \quad \text{and} \quad \frac{\kappa_\beta{}^2}{\left(\kappa_\beta{}^2 + \tau_\beta{}^2\right)^{3/2}}\left(\frac{\tau_\beta}{\kappa_\beta}\right)' = \frac{\kappa^2}{\left(\kappa^2+\tau^2\right)^{3/2}}\left(\frac{\tau}{\kappa}\right)'.$$

From Corollary 4.6, we can give the following theorem without any proof.

**Theorem 4.7.** Let $\alpha$ be a curve in Euclidean 3-space nonzero curvature, $\alpha_T$ be tangent indicatrix of $\alpha$, and $\beta$ be evolute-direction curve of $\alpha_T$. Then

**i)** $\alpha$ is a helix if and only if $\beta$ is a helix.

**ii)** $\alpha$ is a slant helix if and only if $\beta$ is a slant helix.

## 5. Bertrand-direction curves of tangent indicatrix

In this section, we define Bertrand-direction curves of tangent indicatrix and obtain some relations between these curves.

**Definition 5.1.** Let $\alpha$ be a curve in Euclidean 3-space, $\alpha_T$ be tangent indicatrix of $\alpha$, and $\beta$ be $X$-direction curve of $\alpha_T$. If $\beta$ is a Bertrand curve of $\alpha_T$, then $\beta$ is called Bertrand-direction curve of $\alpha_T$. $\alpha_T$ is also called Bertrand-donor curve of $\beta$.

**Proposition 5.2.** Let $\alpha$ be a curve in Euclidean 3-space and $\alpha_T$ be tangent indicatrix of $\alpha$. $\beta$ is Bertrand-direction curve of $\alpha_T$ if and only if the functions in (3.1) are as follows

$$x(s_T) = \cos\theta, \quad y(s_T) = 0, \quad z(s_T) = \sin\theta,$$

where $\theta$ is a constant angle.

**Proof.** From the definition of Bertrand curves, we know that $N_\beta = N_T$. By using equality (3.4), we have the following system of differential equations

$$\left.\begin{aligned} &x' - y\kappa_T = 0 \\ &y' + x\kappa_T - z\tau_T = \kappa_\beta \\ &z' + y\tau_T = 0 \end{aligned}\right\} \tag{5.1}$$

Multiplying the first, second and third equations in (5.1) with $x$, $y$ and $z$, respectively, adding the results, and using equation (3.3), we find $y = 0$. By substituting $y = 0$ into system (5.1), we get $x = c_1$ and $z = c_2$, where $c_1$ and $c_2$ are constant real numbers. Since $X$ is unit, we can give the solution of system (5.1) as follows

$$\{\ x(s_T) = \cos\theta,\ y(s_T) = 0,\ z(s_T) = \sin\theta\ \},$$

where $\theta$ is the constant angle between unit tangent vectors of the curves $\alpha_T$ and $\beta$. □

From Proposition 5.2, we have a method to construct a unit speed Bertrand curve from a unit spherical curve. Moreover, this construction can be achieved just by using Frenet vectors $T_T, B_T$ of tangent indicatrix curve $\alpha_T$ of a given curve $\alpha$, and a constant angle $\theta$ which is the angle between unit tangent vectors of the curves $\alpha_T$ and $\beta$.

Now we can give the following relationships between curvatures of the curve $\alpha_T$ and its Bertrand-direction curve $\beta$.

**Theorem 5.3.** Let $\alpha$ be a curve in Euclidean 3-space, $\alpha_T$ be tangent indicatrix of $\alpha$, and $\beta$ be Bertrand-direction curve of $\alpha_T$. The relations between curvatures of the tangent indicatrix curve $\alpha_T$ and its Bertrand-direction curve $\beta$ can be given as follows

$$\kappa_\beta = \kappa_T \cos\theta - \tau_T \sin\theta \ , \quad \tau_\beta = \kappa_T \sin\theta + \tau_T \cos\theta$$

and

$$\kappa_T = \kappa_\beta \cos\theta + \tau_\beta \sin\theta \ , \quad \tau_T = -\kappa_\beta \sin\theta + \tau_\beta \cos\theta .$$

**Proof.** From the second equation of system (5.1), we can see that $\kappa_\beta = \kappa_T \cos\theta - \tau_T \sin\theta$. From Proposition (5.2), we find $X = T_\beta = \cos\theta T_T + \sin\theta B_T$, and from the definition of Bertrand curves, we have $N_\beta = N_T$. Thus, the unit binormal vector $B_\beta$ of $\beta$ can be found as $B_\beta = T_\beta \times N_\beta = -\sin\theta T_T + \cos\theta B_T$. By differentiating unit binormal vector and using the third equality of (2.3), we find the torsion of $\beta$ as $\tau_\beta = \kappa_T \sin\theta + \tau_T \cos\theta$.

On the other hand, the curvatures of the tangent indicatrix curve $\alpha_T$ can easily be found by using the curvatures of the Bertrand-direction curve $\beta$. □

From Theorem 5.3, we have the following corollary.

**Corollary 5.4.** Let $\alpha$ be a curve in Euclidean 3-space, $\alpha_T$ be tangent indicatrix of $\alpha$, and $\beta$ be Bertrand-direction curve of $\alpha_T$. Then, the following relations between curvatures of the curves $\alpha_T$ and $\beta$ are satisfied

**i)** $$\frac{\tau_T}{\kappa_T} = \frac{-\sin\theta + \frac{\tau_\beta}{\kappa_\beta}\cos\theta}{\cos\theta + \frac{\tau_\beta}{\kappa_\beta}\sin\theta},$$

**ii)** $$\frac{\kappa_T{}^2}{\left(\kappa_T{}^2 + \tau_T{}^2\right)^{3/2}}\left(\frac{\tau_T}{\kappa_T}\right)' = \frac{\kappa_\beta{}^2}{\left(\kappa_\beta{}^2 + \tau_\beta{}^2\right)^{3/2}}\left(\frac{\tau_\beta}{\kappa_\beta}\right)'.$$

From Corollary 5.4, we can give the following theorem without any proof.

**Theorem 5.5.** Let $\alpha$ be a curve in Euclidean 3-space, $\alpha_T$ be tangent indicatrix of $\alpha$, and $\beta$ be Bertrand-direction curve of $\alpha_T$. Then

**i)** $\alpha_T$ is a spherical helix if and only if $\beta$ is a helix.

**ii)** $\alpha_T$ is a spherical slant helix if and only if $\beta$ is a slant helix.

Theorem 5.5 gives a method to construct a helix and a slant helix from a spherical helix and a spherical slant helix, respectively.

Now, by using relations between curvatures of the curve $\alpha$ and its tangent indicatrix $\alpha_T$, we can give relations between curvatures of the curve $\alpha$ and Bertrand-direction curve $\beta$ of tangent indicatrix $\alpha_T$.

**Corollary 5.6.** Let $\alpha$ be a curve in Euclidean 3-space with nonzero curvature, $\alpha_T$ be tangent indicatrix of $\alpha$, and $\beta$ be Bertrand-direction curve of $\alpha_T$. The following relation between curvatures of the curves $\alpha$ and $\beta$ is satisfied

$$\frac{\cos\theta + \left(\tau_\beta/\kappa_\beta\right)\sin\theta}{-\sin\theta + \left(\tau_\beta/\kappa_\beta\right)\cos\theta} = \frac{\kappa^2}{\left(\kappa^2 + \tau^2\right)^{3/2}}\left(\frac{\tau}{\kappa}\right)'.$$

From Corollary 5.6, we can give the following theorem without any proof.

**Theorem 5.7.** Let $\alpha$ be a curve in Euclidean 3-space with nonzero curvature, $\alpha_T$ be tangent indicatrix of $\alpha$, and $\beta$ be Bertrand-direction curve of $\alpha_T$. Then $\alpha$ is a slant helix if and only if $\beta$ is a helix.

## 6. Mannheim-direction curves of tangent indicatrix

In this section, we define Mannheim-direction curves of tangent indicatrix of a curve and obtain some relations between these curves.

**Definition 6.1.** Let $\alpha$ be a curve in Euclidean 3-space, $\alpha_T$ be tangent indicatrix of $\alpha$, and $\beta$ be $X$-direction curve of $\alpha_T$. If $\beta$ is a Mannheim curve of $\alpha_T$, then $\beta$ is called Mannheim-direction curve of $\alpha_T$. $\alpha_T$ is also called Mannheim-donor curve of $\beta$.

**Proposition 6.2.** Let $\alpha$ be a curve in Euclidean 3-space and $\alpha_T$ be tangent indicatrix of $\alpha$. $\beta$ is Mannheim-direction curve of $\alpha_T$ if and only if the functions in (3.1) are as follows

$$x(s_T) = \sin\left(\int \kappa_T ds_T\right), \ y(s_T) = \cos\left(\int \kappa_T ds_T\right), \ z(s_T) = 0.$$

**Proof.** From the definition of Mannheim curves, we know that $N_\beta = B_T$. By using equality (3.4), we have the following system of differential equations

$$\left.\begin{array}{l} x' - y\kappa_T = 0 \\ y' + x\kappa_T - z\tau_T = 0 \\ z' + y\tau_T = \kappa_\beta \end{array}\right\} \quad (6.1)$$

Multiplying the first, second and third equations in (6.1) with $x$, $y$ and $z$, respectively, adding the results, and using equation (3.3), we find $z = 0$. By substituting $z = 0$ into system (6.1), the solution can be obtained as follows

$$\left\{ x(s_T) = \sin\left(\int \kappa_T ds_T\right), y(s_T) = \cos\left(\int \kappa_T ds_T\right), z(s_T) = 0 \right\}. \ \square$$

From Proposition 6.2, we have a method to construct a unit speed Mannheim curve from a unit spherical curve. Moreover, this construction can be achieved just by using unit Frenet vectors $T_T$ and $N_T$, and curvature $\kappa_T$ of tangent indicatrix curve $\alpha_T$ of a given curve $\alpha$.

Now we can give the following relationships between curvatures of the curve $\alpha_T$ and its Mannheim-direction curve $\beta$.

**Theorem 6.3.** Let $\alpha$ be a curve in Euclidean 3-space, $\alpha_T$ be tangent indicatrix of $\alpha$, and $\beta$ be Mannheim-direction curve of $\alpha_T$. The relations between curvatures of the tangent indicatrix curve $\alpha_T$ and its Mannheim-direction curve $\beta$ can be given as follows

$$\kappa_\beta = \tau_T \cos\left(\int \kappa_T ds_T\right), \quad \tau_\beta = \tau_T \sin\left(\int \kappa_T ds_T\right)$$

and

$$\kappa_T = \frac{\kappa_\beta^{\ 2}}{\kappa_\beta^{\ 2} + \tau_\beta^{\ 2}}\left(\frac{\tau_\beta}{\kappa_\beta}\right)', \quad \tau_T = \sqrt{\kappa_\beta^{\ 2} + \tau_\beta^{\ 2}}.$$

**Proof.** From the third equation of system (6.1), we can easily see that $\kappa_\beta = \tau_T \cos\left(\int \kappa_T ds_T\right)$. From Proposition (6.2), we find $X = T_\beta = \sin\left(\int \kappa_T ds_T\right)T_T + \cos\left(\int \kappa_T ds_T\right)N_T$, and from the definition of Mannheim curves, we have $N_\beta = B_T$. The unit binormal vector $B_\beta$ of the curve $\beta$ can be found as $B_\beta = T_\beta \times N_\beta = \cos\left(\int \kappa_T ds_T\right)T_T - \sin\left(\int \kappa_T ds_T\right)N_T$. By differentiating unit binormal vector and using the third equality of (2.3), we get $\tau_\beta = \tau_T \sin\left(\int \kappa_T ds_T\right)$.

On the other hand, by a simple computation with curvatures in Theorem 6.3, we can see the torsion of Mannheim-donor curve $\alpha_T$ as $\tau_T = \sqrt{\kappa_\beta^2 + \tau_\beta^2}$. Moreover, from Theorem 6.3, we have $\frac{\tau_\beta}{\kappa_\beta} = \tan\left(\int \kappa_T ds_T\right)$, and so $\int \kappa_T ds_T = \arctan\left(\frac{\tau_\beta}{\kappa_\beta}\right)$. By differentiating the last equality, we get $\kappa_T = \frac{\kappa_\beta^2}{\kappa_\beta^2 + \tau_\beta^2}\left(\frac{\tau_\beta}{\kappa_\beta}\right)'$. □

From Theorem 6.3, we have the following corollary.

**Corollary 6.4.** Let $\alpha$ be a curve in Euclidean 3-space, $\alpha_T$ be tangent indicatrix of $\alpha$, and $\beta$ be Mannheim-direction curve of $\alpha_T$. Then the following relation between curvatures of the curves $\alpha_T$ and $\beta$ is satisfied

$$\frac{\tau_T}{\kappa_T} = \mp \frac{1}{\dfrac{\kappa_\beta^2}{\left(\kappa_\beta^2 + \tau_\beta^2\right)^{3/2}}\left(\dfrac{\tau_\beta}{\kappa_\beta}\right)'} .$$

From Theorem 6.3 and Corollary 6.4, we can give the following theorem without any proof.

**Theorem 6.5** Let $\alpha$ be a curve in Euclidean 3-space, $\alpha_T$ be tangent indicatrix of $\alpha$, and $\beta$ be Mannheim-direction curve of $\alpha_T$. Then

**i)** $\alpha_T$ is a spherical helix if and only if $\beta$ is a slant helix.

**ii)** $\alpha_T$ is a circle or a part of circle on unit sphere if and only if $\beta$ is a straight line.

Theorem 6.5 (i) gives a method to construct a slant helix from a spherical helix.

Now, by using relations between curvatures of the curve $\alpha$ and its tangent indicatrix $\alpha_T$, we can give the following relations between curvatures of the curve $\alpha$ and Mannheim-direction curve $\beta$ of the tangent indicatrix $\alpha_T$ of $\alpha$.

**Corollary 6.6.** Let $\alpha$ be a curve in Euclidean 3-space, $\alpha_T$ be tangent indicatrix of $\alpha$, and $\beta$ be Mannheim-direction curve of $\alpha_T$. The following relations between curvatures of the curves $\alpha$ and $\beta$ are satisfied

$$\frac{\tau_\beta}{\kappa_\beta} = \tan\left(\int \sqrt{\kappa^2 + \tau^2} ds\right) \quad \text{and} \quad \frac{\kappa_\beta^2}{\left(\kappa_\beta^2 + \tau_\beta^2\right)^{3/2}}\left(\frac{\tau_\beta}{\kappa_\beta}\right)' = \mp \frac{1}{\dfrac{\kappa^2}{\left(\kappa^2 + \tau^2\right)^{3/2}}\left(\dfrac{\tau}{\kappa}\right)'} .$$

From Corollary 6.6, we can give the following theorem without any proof.

**Theorem 6.7.** Let $\alpha$ be a curve in Euclidean 3-space, $\alpha_T$ be tangent indicatrix of $\alpha$, and $\beta$ be Mannheim-direction curve of $\alpha_T$. Then **i)** $\alpha$ is a straight line if and only if $\beta$ is a helix.
**ii)** $\alpha$ is a slant helix if and only if $\beta$ is a slant helix.

## 7. Examples

In this section, we give two examples of direction curves of tangent indicatrices. In both examples, we find evolute-direction, Bertrand-direction, and Mannheim-direction curves of tangent indicatrix of a given curve and illustrate these curves by using the program Wolfram Mathematica 9.0.

**Example 7.1.** Let consider a space curve given by the parametrization $\alpha(s)=\left(\cos\frac{s}{\sqrt{2}},\sin\frac{s}{\sqrt{2}},\frac{s}{\sqrt{2}}\right)$, where $s$ is the arc-length parameter of $\alpha$. The tangent indicatrix of the curve $\alpha$ can be found as $\alpha_T(s_T)=\left(-\frac{1}{\sqrt{2}}\sin\left(\sqrt{2}s_T\right),\frac{1}{\sqrt{2}}\cos\left(\sqrt{2}s_T\right),\frac{1}{\sqrt{2}}\right)$, where $s_T$ is the arc-length parameter of the tangent indicatrix curve $\alpha_T$. The evolute-direction, Bertrand-direction, and Mannheim-direction curves of the tangent indicatrix curve $\alpha_T$ are obtained, respectively, as follows

$$\beta_E(s_T)=\left(-\frac{1}{\sqrt{2}}\sin\theta_1\cos\left(\sqrt{2}s_T\right)+c_1,-\frac{1}{\sqrt{2}}\sin\theta_1\sin\left(\sqrt{2}s_T\right)+c_2,s_T\cos\theta_1+c_3\right),$$

$$\beta_B(s)=\left(-\frac{1}{\sqrt{2}}\cos\theta\sin\left(\sqrt{2}s_T\right)+c_4,\frac{1}{\sqrt{2}}\cos\theta\cos\left(\sqrt{2}s_T\right)+c_5,s_T\sin\theta+c_6\right),$$

$$\beta_M(s)=\left(-s_T\sin\theta_2+c_7,-s_T\cos\theta_2+c_8,c_9\right),$$

where $\theta_1,\theta_2,c_i,(i=1,...,9)$ are real integration constants and $\theta$ is a real constant angle. By taking $\theta_1=\theta_2=\pi/4$, $c_i=0$, $(i=1,...,9)$ and $\theta=\pi/3$, we can illustrate the curve $\alpha$, the tangent indicatrix curve $\alpha_T$, and evolute-direction, Bertrand direction, and Mannheim-direction curves of the curve $\alpha_T$ in Figures 1, 2, 3, 4, 5, respectively.

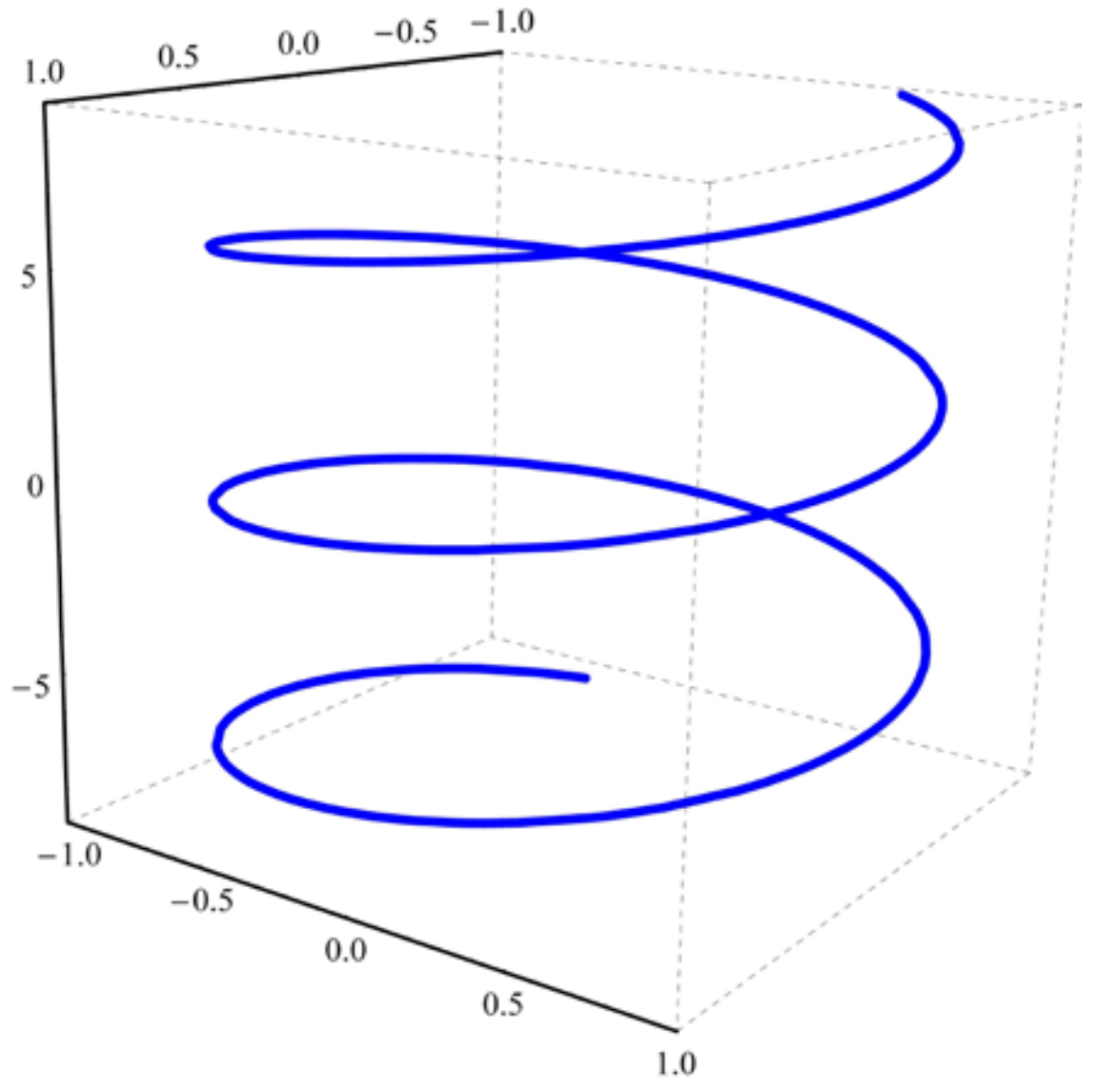

**Figure 1** Given curve $\alpha$

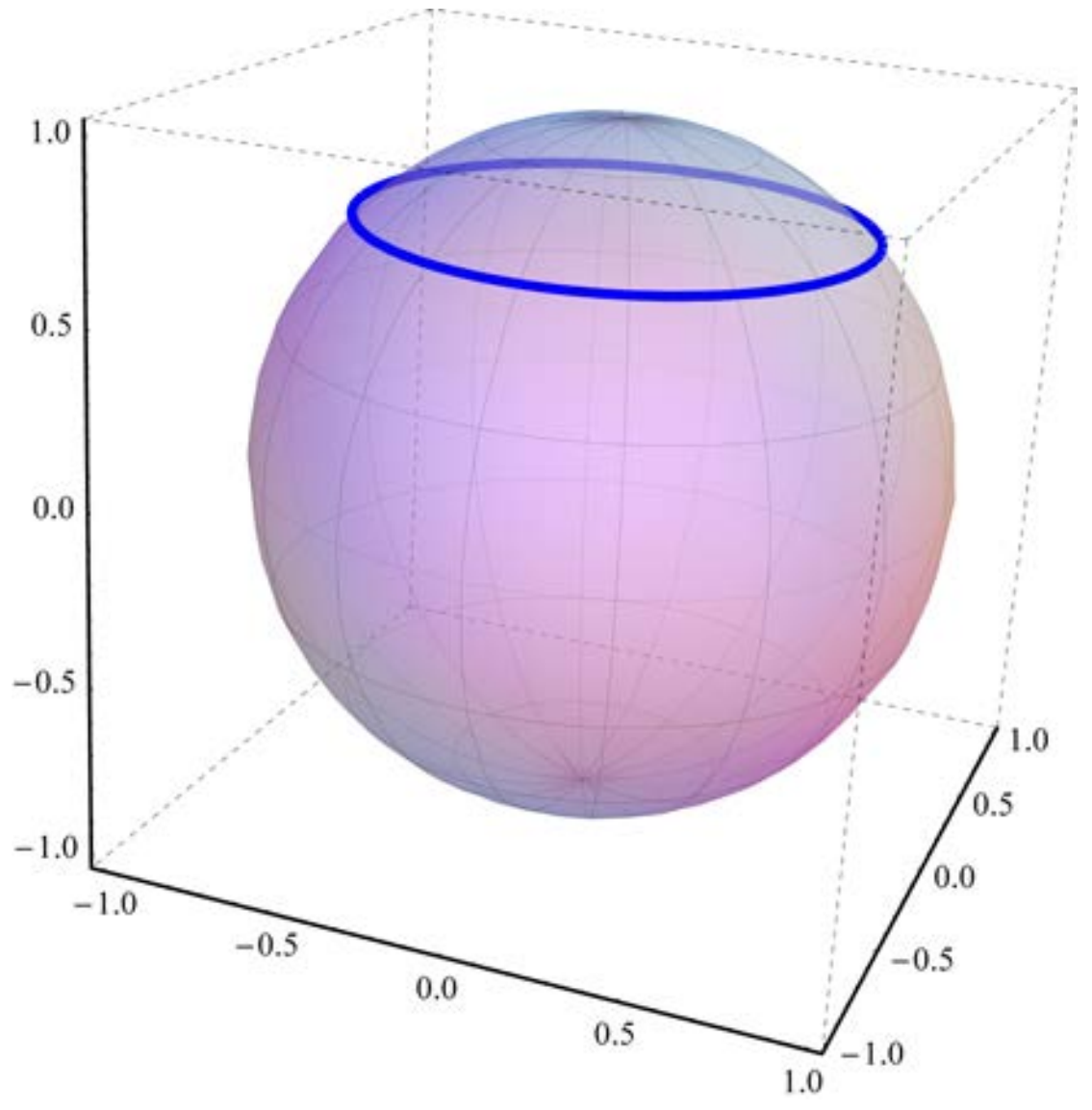

**Figure 2** Tangent indicatrix of the curve $\alpha$

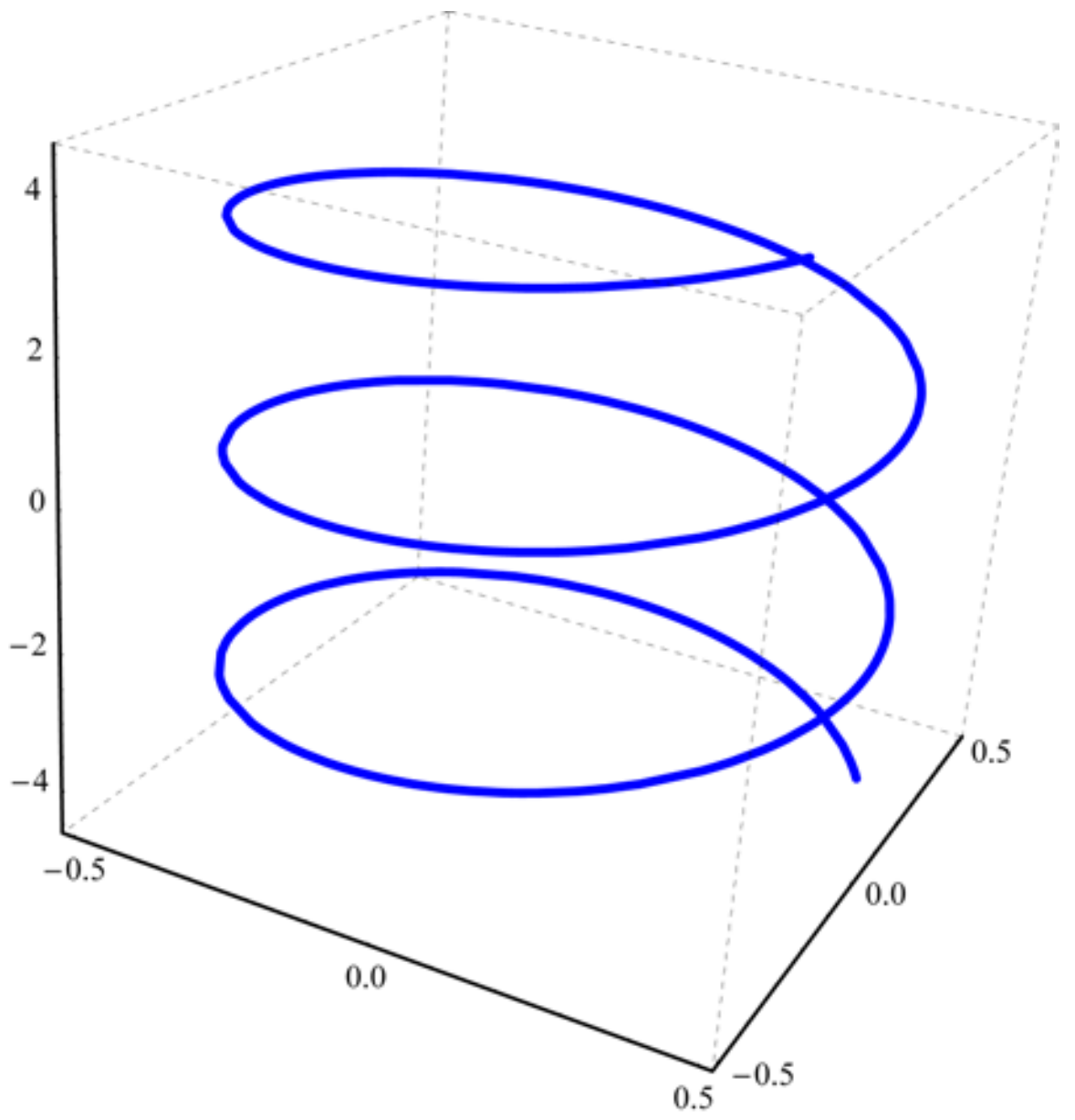

**Figure 3** Evolute-direction curve of the tangent Indicatrix

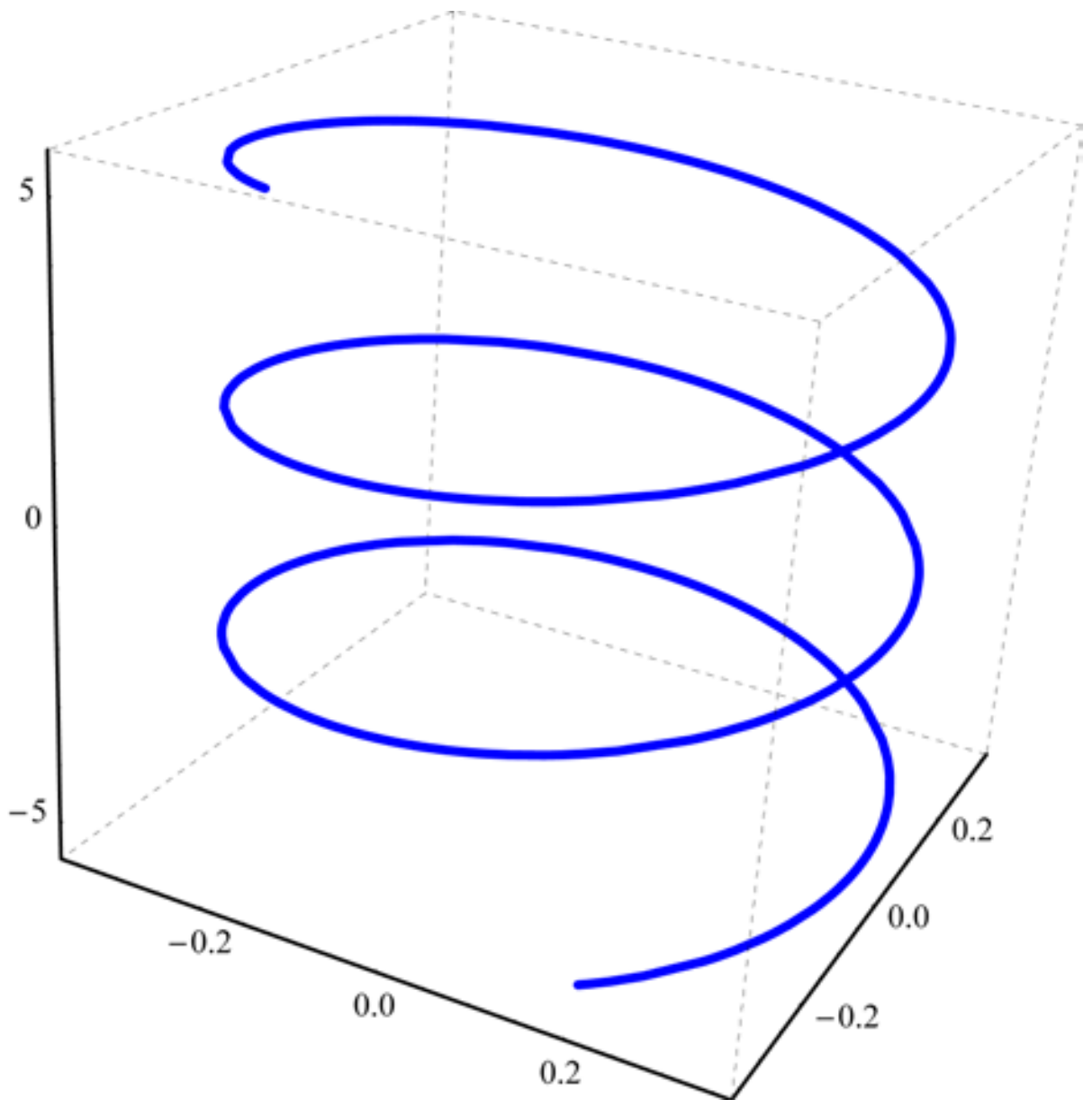

**Figure 4** Bertrand-direction curve of the tangent indicatrix

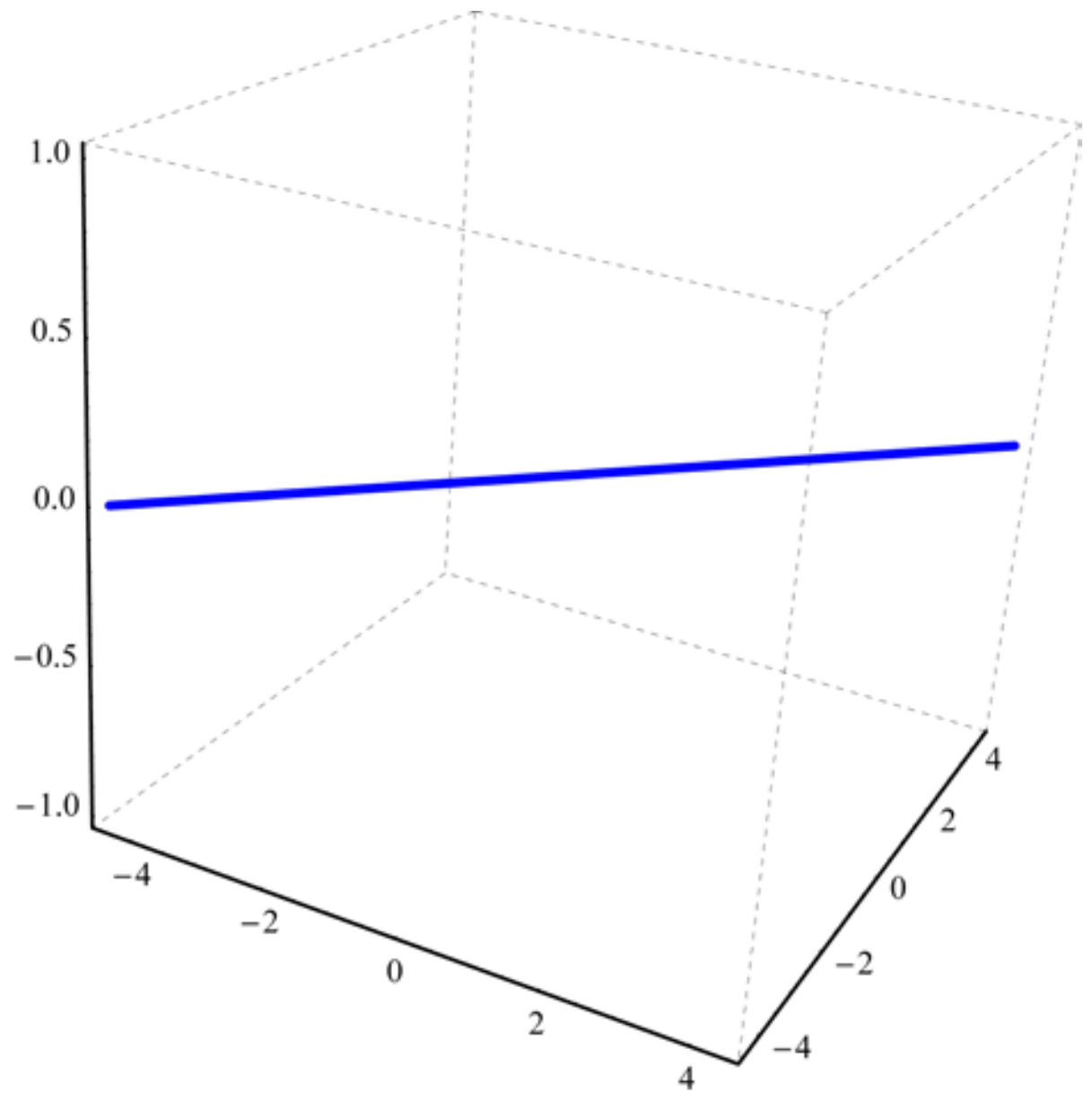

**Figure 5** Mannheim-direction curve of the tangent indicatrix

**Example 7.2.** Let consider a space curve given by the parametrization

$$\alpha(t)=\left(\frac{3}{\sqrt{2}}\sin\left(\sqrt{2}\,t\right)\cos t-2\sin t\cos\left(\sqrt{2}\,t\right),\frac{3}{\sqrt{2}}\cos\left(\sqrt{2}\,t\right)\cos t+2\sin t\sin\left(\sqrt{2}\,t\right),-\frac{1}{\sqrt{2}}\cos t\right).$$

The tangent indicatrix of the curve $\alpha$ can be found as

$$\alpha(t)=\left(\cos\left(\sqrt{2}\,t\right)\cos t+\frac{1}{\sqrt{2}}\sin\left(\sqrt{2}\,t\right)\sin t,-\sin\left(\sqrt{2}\,t\right)\cos t+\frac{1}{\sqrt{2}}\cos\left(\sqrt{2}\,t\right)\sin t,\frac{1}{\sqrt{2}}\sin t\right).$$

The evolute-direction, Bertrand-direction, and Mannheim-direction curves of the tangent indicatrix curve $\alpha_T$ can be found as

$$\beta_E(t)=\Bigg(\frac{1}{8}\Big(-\big(-1+\sqrt{2}\big)\cos\big(\theta_1-\sqrt{2}\,t\big)+\big(1+\sqrt{2}\big)\cos\big(\theta_1+\sqrt{2}\,t\big)$$
$$-\big(3+2\sqrt{2}\big)\cos\Big(\theta_1+\big(-2+\sqrt{2}\big)t\Big)+\big(-3+2\sqrt{2}\big)\cos\Big(\theta_1-\big(2+\sqrt{2}\big)t\Big)\Big)+c_1,$$
$$\frac{1}{8}\Big(-\big(-1+\sqrt{2}\big)\sin\big(\theta_1-\sqrt{2}\,t\big)-\big(1+\sqrt{2}\big)\sin\big(\theta_1+\sqrt{2}\,t\big)$$
$$+\big(3+2\sqrt{2}\big)\sin\Big(\theta_1+\big(-2+\sqrt{2}\big)t\Big)+\big(-3+2\sqrt{2}\big)\sin\Big(\theta_1-\big(2+\sqrt{2}\big)t\Big)\Big)+c_2,$$
$$\frac{1}{4\sqrt{2}}\big(-2t\cos(\theta_1)+\sin(\theta_1-2t)\big)+c_3\Bigg),$$
$$\beta_B(t)=\Bigg(\frac{1}{2}\Big(2\cos t\cos\big(\sqrt{2}\,t\big)+\sqrt{2}\sin t\sin\big(\sqrt{2}\,t\big)\Big)\big(\cos\theta+\sin\theta\big)+c_4,$$
$$\frac{1}{2}\Big(\sqrt{2}\cos\big(\sqrt{2}\,t\big)\sin t-2\cos t\sin\big(\sqrt{2}\,t\big)\Big)\big(\cos\theta+\sin\theta\big)+c_5,$$
$$\frac{1}{\sqrt{2}}\sin t\big(\cos\theta-\sin\theta\big)+c_6\Bigg),$$
$$\beta_M(t)=\Bigg(\frac{1}{8}\Big(\big(1+\sqrt{2}\big)\sin\big(\theta_2-\sqrt{2}\,t\big)-\big(-1+\sqrt{2}\big)\sin\big(\theta_2+\sqrt{2}\,t\big)$$
$$-\big(3+2\sqrt{2}\big)\sin\Big(\theta_2+\big(2-\sqrt{2}\big)t\Big)+\big(-3+2\sqrt{2}\big)\sin\Big(\theta_2+\big(2+\sqrt{2}\big)t\Big)\Big)+c_7,$$
$$\frac{1}{8}\Big(-\big(1+\sqrt{2}\big)\cos\big(\theta_2-\sqrt{2}\,t\big)-\big(-1+\sqrt{2}\big)\cos\big(\theta_2+\sqrt{2}\,t\big)$$
$$+\big(3+2\sqrt{2}\big)\cos\Big(\theta_2+\big(2-\sqrt{2}\big)t\Big)+\big(-3+2\sqrt{2}\big)\cos\Big(\theta_2+\big(2+\sqrt{2}\big)t\Big)\Big)+c_8,$$
$$-\frac{1}{4\sqrt{2}}\big(\cos(\theta_2+2t)-2t\sin(\theta_2)\big)+c_9\Bigg),$$

respectively, where $c_i,(i=1,...,9)$ and $\theta_1,\theta_2$ are real integration constants and $\theta$ is a real constant angle. By taking $c_i=0,(i=1,...,9)$, $\theta=\pi/3$, $\theta_1=\pi/3$, $\theta_2=0$, we can illustrate the curve $\alpha$, the tangent indicatrix curve $\alpha_T$ of $\alpha$, and its evolute-direction, Bertrand direction, and Mannheim-direction curves in Figures 6, 7, 8, 9, 10, respectively.

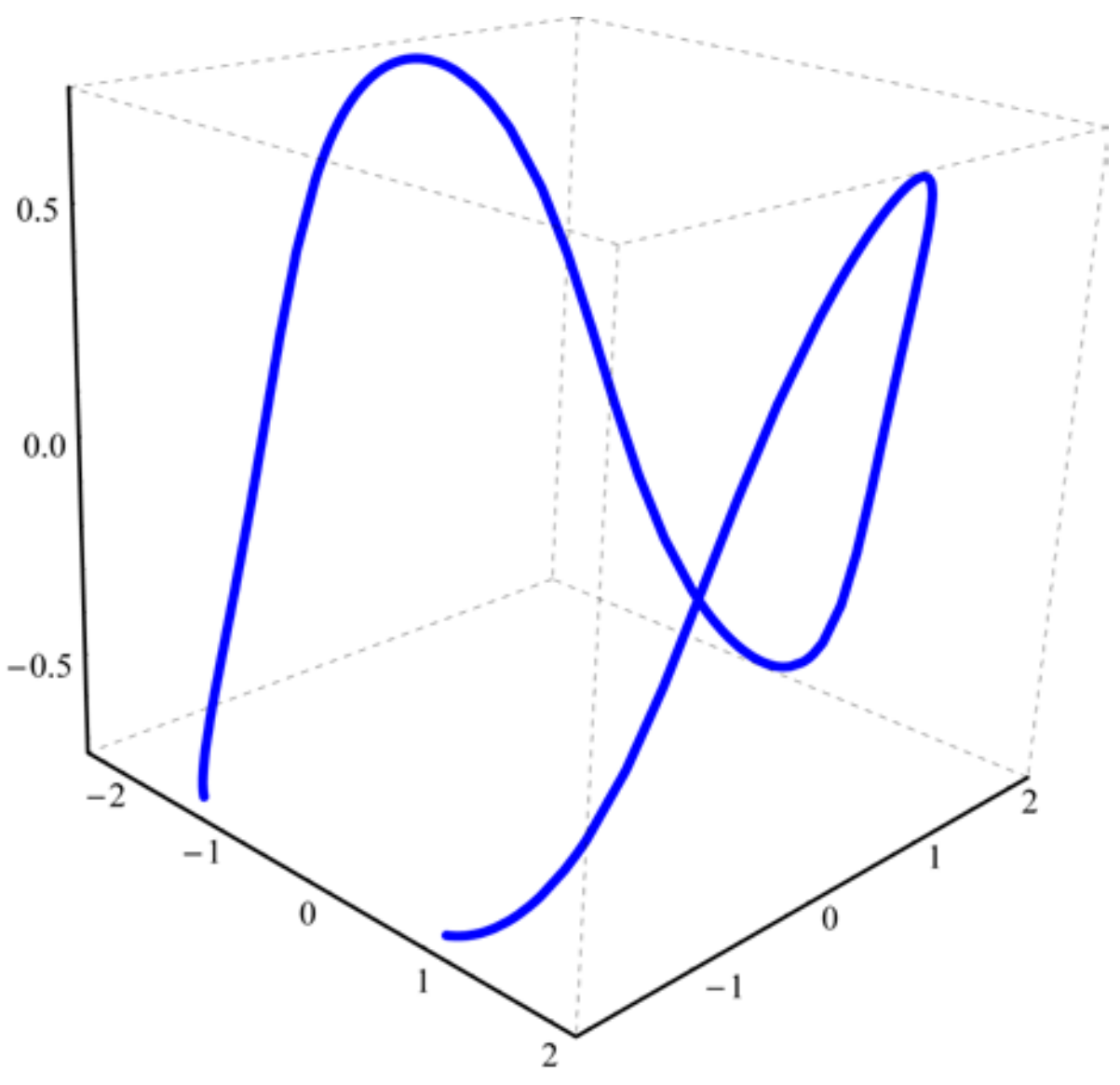

**Figure 6** Given curve $\alpha$

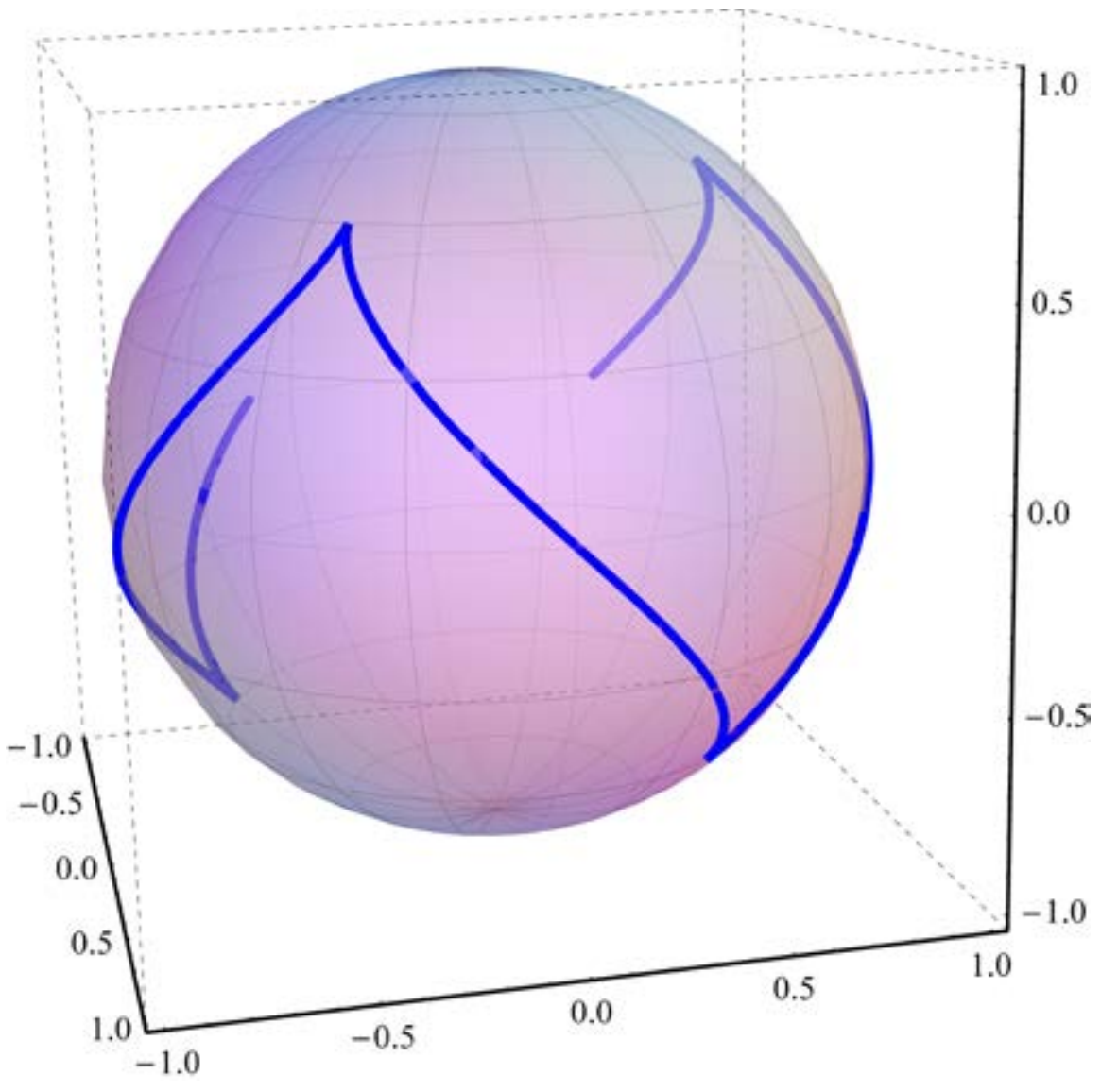

**Figure 7** Tangent indicatrix of the curve $\alpha$

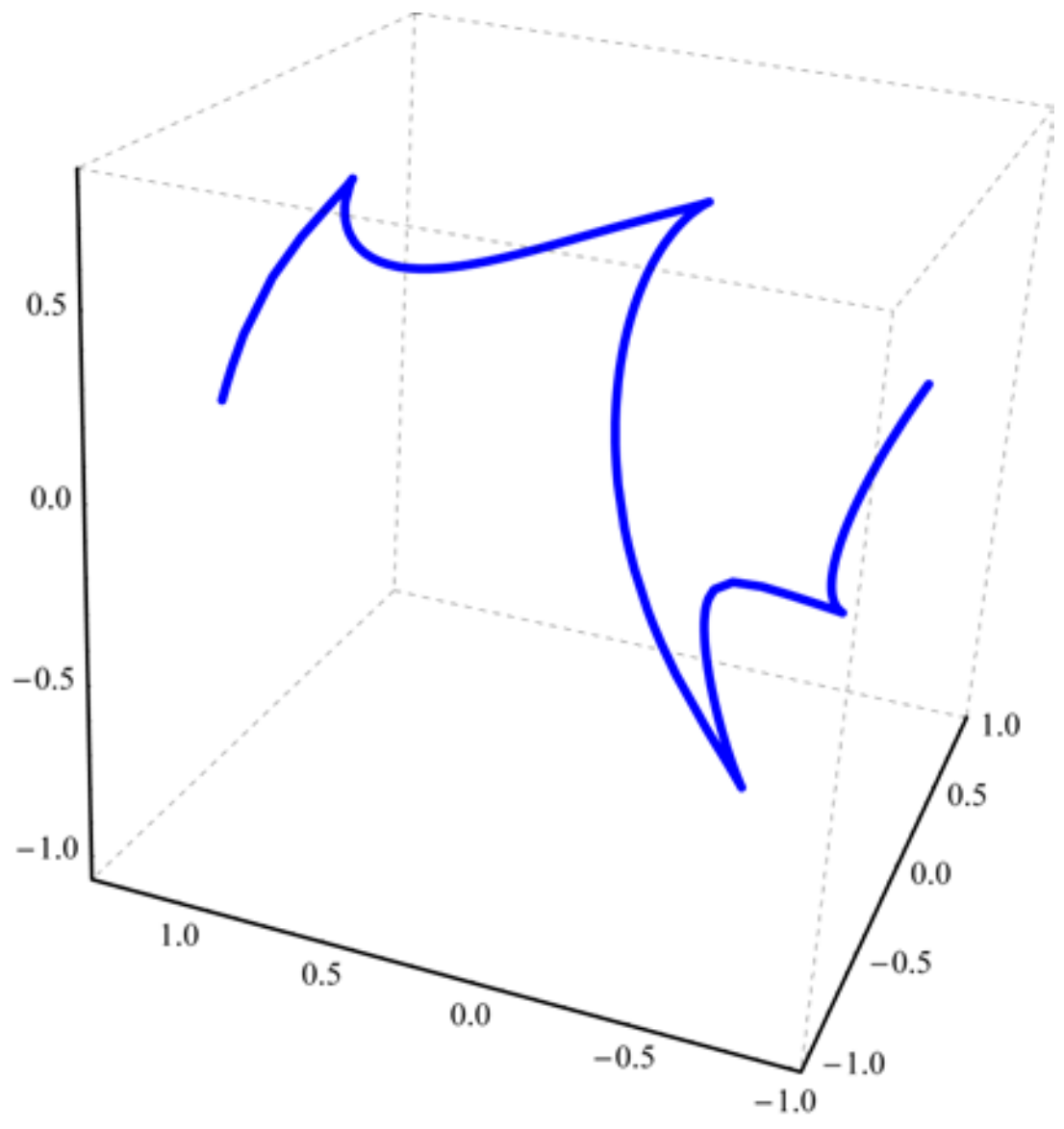

**Figure 8** Evolute-direction curve of the tangent indicatrix

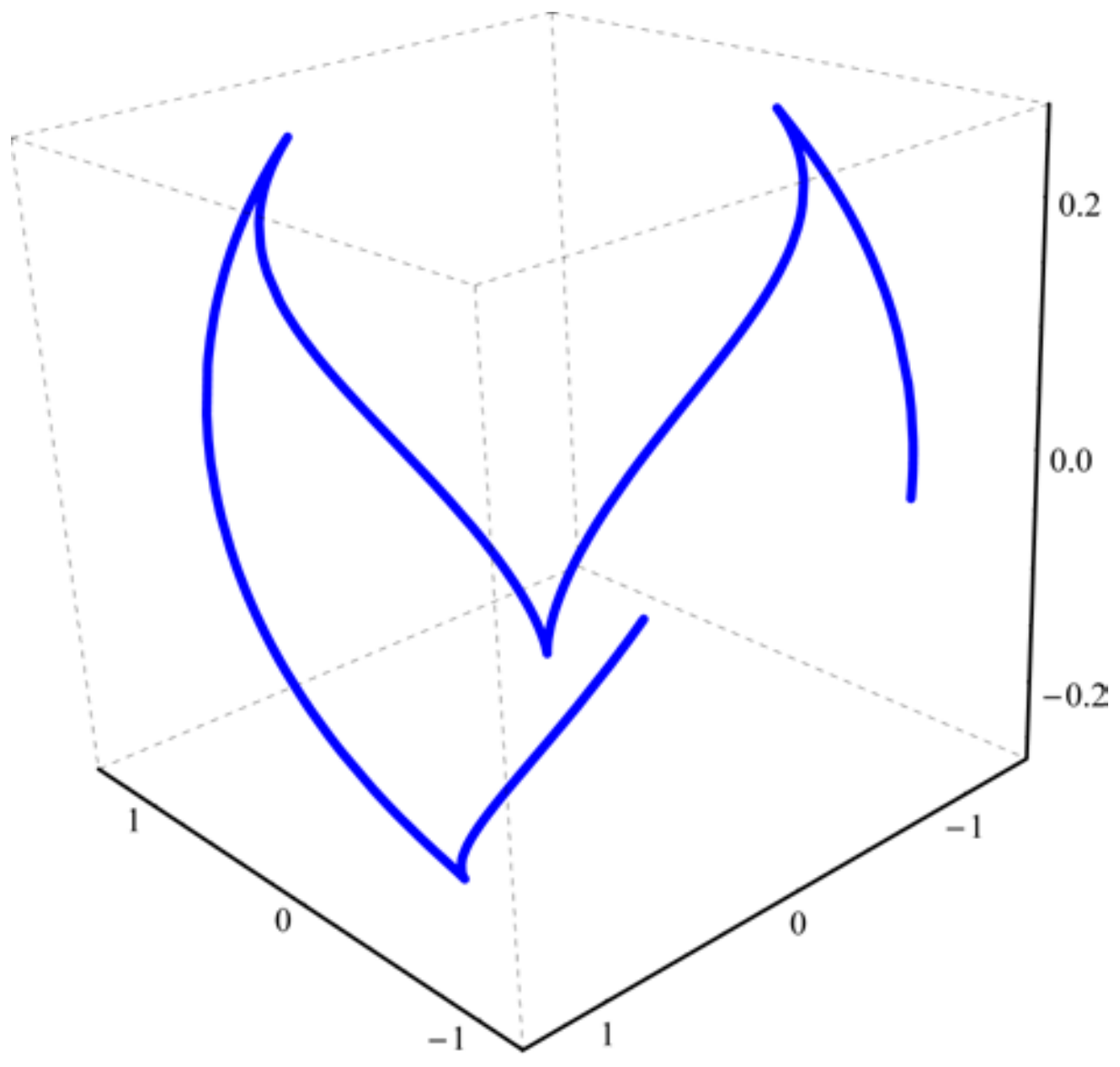

**Figure 9** Bertrand-direction curve of the tangent indicatrix

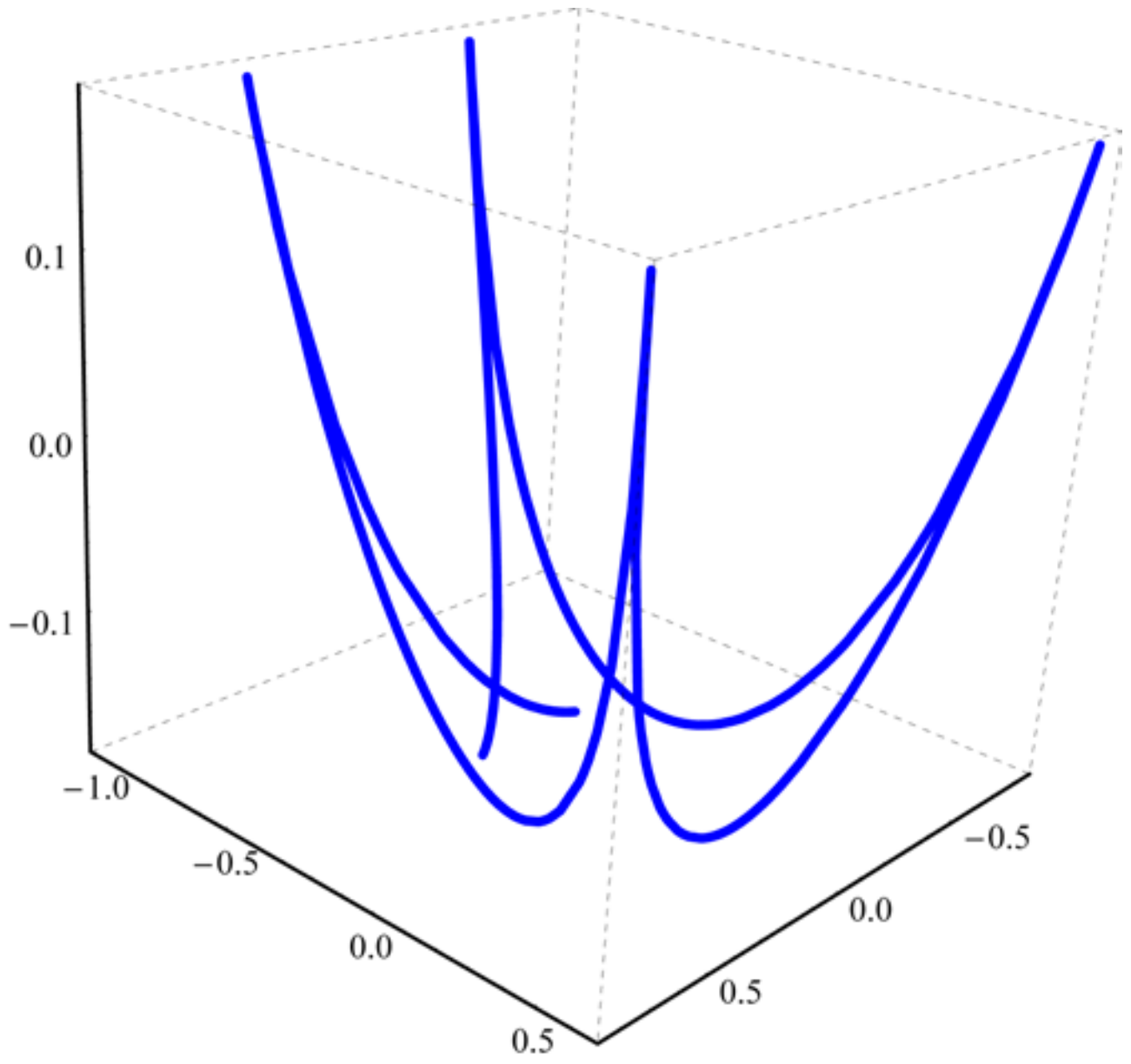

**Figure 10** Mannheim-direction curve of the tangent indicatrix